\documentclass[12pt,reqno]{article}

\usepackage[usenames]{color}
\usepackage{amssymb}
\usepackage{amsmath}
\usepackage{amsthm}
\usepackage{amsfonts}
\usepackage{amscd}
\usepackage{graphicx}

\usepackage[colorlinks=true,
linkcolor=webgreen,
filecolor=webbrown,
citecolor=webgreen]{hyperref}

\definecolor{webgreen}{rgb}{0,.5,0}
\definecolor{webbrown}{rgb}{.6,0,0}

\usepackage{color}

\usepackage{graphics}
\usepackage{latexsym}

\begin{document}

\theoremstyle{plain}
\newtheorem{theorem}{Theorem}
\newtheorem{corollary}[theorem]{Corollary}
\newtheorem{lemma}[theorem]{Lemma}
\newtheorem{proposition}[theorem]{Proposition}

\theoremstyle{definition}
\newtheorem{definition}[theorem]{Definition}
\newtheorem{example}[theorem]{Example}
\newtheorem{conjecture}[theorem]{Conjecture}

\theoremstyle{remark}
\newtheorem{remark}[theorem]{Remark}

\begin{center}
\vskip 1cm{\LARGE\bf  On a New Congruence in the Catalan Triangle
\vskip 1cm}
\large
Jovan Miki\'{c}\\
Faculty of Technology\\
Faculty of Natural Sciences and Mathematics\\
University of Banja Luka\\
Bosnia and Herzegovina\\
\href{mailto:jovan.mikic@tf.unibl.org}{\tt jovan.mikic@tf.unibl.org} \\
\end{center}

\vskip .2in

\begin{abstract}
For $0\leq k \leq n$, the number $C(n,k)$ represents the number of all lattice paths
in the plane from the point $(0,0)$ to the point $(n,k)$, using steps $(1,0)$ and $(0,1)$, that never rise above the main diagonal $y=x$. The Fuss-Catalan number of order three $C^{(3)}_n$ represents the number of all lattice paths in the plane from the point $(0,0)$
to the point $(2n,n)$, using steps  $(1,0)$ and $(0,1)$, that do not rise above the line $y=\frac{x}{2}$. We present a new alternating convolution formula
for the numbers $C(2n,k)$. By using a new class of binomial sums that we call $M$ sums, we prove that this sum is divisible by $C^{(3)}_n$ and by the central binomial coefficient $\binom{2n}{n}$. We do this by examining the numbers $T(n,j)=\frac{1}{2n+1}\binom{2n+j}{j}\binom{2n+1}{n+j+1}$, for which we present a new combinatorial interpretation, connecting them to the generalized Schr\"{o}der numbers of order two.
\end{abstract}

\section{Introduction}\label{s:1}

 Let $C_n=\frac{1}{n+1}\binom{2n}{n}$ denote the $n$-th Catalan number. Catalan numbers form the famous sequence \cite{RS} with the most applications in combinatorics, after the binomial coefficients. For example, $C_n$ is the number of all lattice paths in the plane from the point $(0,0)$ to the point $(n,n)$, using steps $(1,0)$ and $(0,1)$, that never rise above the main diagonal $y=x$.

 This sequence starts with: $C_0=1$, $C_1=1$, $C_2=2$, $C_3=5$, $C_4=14$, etc., and it can be found as sequence $A000108$ in \cite{NS}. See also \cite[Introduction]{JCA}.

 Let $n$ and $k$ be non-negative integers such that $k\leq n$. The number $C(n,k)$ represents the number of all lattice paths in the plane from the point $(0,0)$
 to the point $(n,k)$, using steps $(1,0)$ and $(0,1)$, that never rise above the main diagonal $y=x$. The numbers $C(n,k)$ therefore form the Catalan triangle, where the Catalan numbers appear on the main diagonal:
 \[
 \begin{matrix}
 1\\
 1 & 1\\
 1& 2 &2\\
 1&3&5&5\\
 1& 4 &9 &14 &14\\
 1&5&14& 28& 42 &42\\
 1&6&20&48&90&132&132\\
 \ldots
 \end{matrix}
 \]

 Furthermore, the following two recurrences hold for $C(n,k)$:
 \begin{align*}
     C(n+1,k)&=C(n+1,k-1)+C(n,k)\text{, } 1<k <n+1\text{  and}\\
     C(n+1,n+1)&=C(n+1,n)\text{.}
 \end{align*}
 
 It is well-known that $\sum_{k=0}^{n}C(n,k)=C_{n+1}$. It follows that the sum of elements of every row of the Catalan triangle is again Catalan number.

 An explicit expression for $C(n,k)$ is
 \begin{equation}\label{eq:1}
     C(n,k)=\frac{n-k+1}{n+1}\binom{n+k}{n}\text{.}
 \end{equation}

 The number $C(n,k)$ is always an integer due to the fact that
 \begin{equation}\label{eq:2}
     C(n,k)=\binom{n+k}{n}-\binom{n+k}{n+1}\text{.}
 \end{equation}

 The numbers $C(n,k)$ \cite[Section 15, p.\ 347]{TC} can also be defined in terms of binary words.
 
 There exist other identities with the numbers $C(n,k)$. See for example \cite[Thms.~1.1,1.2; p.\ 2]{LO}.
 Note that there are some other triangles which are also called Catalan triangles. For example, there is an another Catalan triangle \cite{LS} introduced by Shapiro. See also \cite[Section 14, p.\ 333]{TC} and \cite{MOR}.

 Let $C^{(3)}_n=\frac{1}{2n+1}\binom{3n}{n}$ denote the Fuss-Catalan number of order three. The number $C^{(3)}_n$ counts all lattice paths from the point $(0,0)$ to the point $(2n,n)$, using steps $(1,0)$ and $(0,1)$, that do not rise above the line $y=\frac{x}{2}$. See also \cite[Eq.~(2.2), p.\ 5]{JCA}. This sequence starts with: $C^{(3)}_0=1$, $C^{(3)}_1=1$, $C^{(3)}_2=3$, $C^{(3)}_3=12$, $C^{(3)}_4=55$, etc., and it can be found as sequence $A001764$ in \cite{NS}.

  The central binomial coefficient $\binom{2n}{n}$ represents the number of all lattice paths in the plane from the point $(0,0)$ to the point $(n,n)$, using steps $(1,0)$ and $(0,1)$. This sequence starts with $1$, $2$, $6$, $20$, $70$, etc., and it can be found as sequence $A000984$ in \cite{NS}.

  In this article, we present a connection between the numbers $C(2n,k)$ and $C^{(3)}_n$. Let $n$ be a non-negative integer, and let $m$ be a natural number. Let us consider the following alternating convolution for the numbers $C(2n,k)$:
  \begin{equation}\label{eq:3}
      S(2n,m)=\sum_{k=0}^{2n}(-1)^k\binom{2n}{k}^m C(2n,k)\cdot C(2n,2n-k)\text{.}
  \end{equation}

  We present a new alternating convolution formula in the Catalan triangle as our first main result.
  \begin{theorem}\label{t:1}
  For any $n \in \mathbb{N}_0$, we have
  \begin{equation}\label{eq:4}
      S(2n,1)=(-1)^n\cdot C^{(3)}_n\cdot C_n\cdot (2n^2+n+1)\text{.}
  \end{equation}
  \end{theorem}

  Our other two main results are the following theorems.

  \begin{theorem}\label{t:2}
  The sum $S(2n,m)$ is divisible by $ C^{(3)}_n$ for any non-negative integer $n$ and for any natural number $m$.
  \end{theorem}

  \begin{theorem}\label{t:3}
  The sum $S(2n,m)$ is divisible by $\binom{2n}{n}$ for any non-negative integer $n$ and for any natural number $m$.
  \end{theorem}

  For proofs of Thms.~\ref{t:1}, \ref{t:2},  and \ref{t:3} we use a new class of binomial sums that we call $M$ sums that include numbers $T(n,j)=\frac{1}{2n+1}\binom{2n+j}{j}\binom{2n+1}{n+j+1}$. The numbers $T(n,j)$ are integers for all non-negative integers $j$ such that $j \leq n$. Let us recall that, for $m \geq 2n$, the generalized Schr\"{o}der numbers $Schr(n,m,2)$ of order two represent the number of all lattice paths in the plane from the point $(0,0)$ to the point $(n,m)$, using steps $(1,0)$, $(0,1)$, and $(1,1)$,
  that never go below the line $y=2x$. Furthermore, let $Schr(n,m,j,2)$ denote the number of all lattice paths
in the plane from the point $(0,0)$ to the point $(n,m)$, using steps $(1,0)$, $(0,1)$, and $(1,1)$,
that never go below the line $y=2x$, and with exactly $j$ steps of the form $(1,0)$.

Our fourth main result is a new combinatorial interpretation of $T(n,j)$.

\begin{theorem}\label{t:4}
Let $n$ be a natural number, and let $j$ be a non-negative integer such that $j\leq n$.
The number $T(n,j)=\frac{1}{2n+1}\binom{2n+j}{j}\binom{2n+1}{n+j+1}$ represents the number of all lattice paths
in the plane from the point $(0,0)$ to the point $(n,2n)$, using steps $(1,0)$, $(0,1)$, and $(1,1)$,
that never go below the line $y=2x$ and with exactly $j$ steps of the form $(1,0)$. 
\end{theorem}

Finally, we shall prove a more general result than Theorem \ref{t:4}.
\begin{theorem}\label{t:5}
Let $n$ be a natural number, and let $j$ be a non-negative integer such that $j\leq n$. Let $m$ be a natural number such that $m\geq 2n$. Then
\begin{equation}\label{eq:5}
Schr(n,m,j,2)=\frac{m-2n+1}{n}\binom{n}{j}\binom{m+j}{n-1}\text{.}
\end{equation}
\end{theorem}

\section{A New Class of Binomial Sums}\label{s:2}

In this section we define a crucial tool in our research, namely the $M$ sums.
We start with the following definition.
\begin{definition}\label{def:1}
Let $S(n,m,a)=\sum_{k=0}^{n}\binom{n}{k}^m\cdot F(n,k,a)$, where $F(n,k,a)$ is an integer-valued function and the number $a$ is also an integer.
Then $M$ sums for the sum $S(n,m,a)$ are given as follows:
\begin{equation}\label{eq:6}
 M_S(n,j,t;a)=\binom{n-j}{j}\sum_{v=0}^{n-2j}\binom{n-2j}{v}\binom{n}{j+v}^t F(n,j+v,a)\text{,}\\
\end{equation}
where $j$ and $t$ are non-negative integers such that $j\leq \lfloor\frac{n}{2}\rfloor$.
\end{definition}

Obviously, 
\begin{equation}\label{eq:7}
    S(n,m,a)=M_S(n,0,m-1;a)\text{.}
\end{equation}

It is known that the $M$ sums \cite[Eq.~(\ref{eq:8}), p.\ 3]{JM1} satisfy the following main recurrence:
\begin{equation}\label{eq:8}
M_S(n,j,t+1;a)=\binom{n}{j}\sum_{u=0}^{\lfloor \frac{n-2j}{2}\rfloor}\binom{n-j}{u} M_S(n,j+u,t;a)\text{.}
\end{equation}

Let $S(n,m,a)=\sum_{k=0}^{n}\binom{n}{k}^m F(n,k,a)$ be a sum from Definition \ref{def:1}, and let $P(n,m,a)=\sum_{k=0}^{n}\binom{n}{k}^m\binom{a+k}{a}\binom{a+n-k}{a} F(n,k,a)$.

Then, the following equation is true \cite[Thm.\ 3, p.\ 3]{JM1}:
\begin{equation}\label{eq:9}
M_P(n,j,0;a)=\binom{a+j}{a}\sum_{l=0}^{a}\binom{n-j+l}{l}\binom{n-j}{a-l}M_S(n,j+a-l,0;a)\text{.}
\end{equation}

We use the following agreement.
If the sum $S(n,m,a)$ does not depend on $a$, we shall write $S(n,m)$ instead of $S(n,m,a)$, and $F(n,k)$ instead of $F(n,k,a)$. Similarly, in that case, we shall write $M_S(n,j,0)$ instead of $M_S(n,j,0;a)$.

Let $S(2n,m)$ be our main sum from Eq.~(\ref{eq:3}).

We shall prove the following lemma.
\begin{lemma}\label{l:1}
\begin{equation}\label{eq:10}
M_S(2n,j,0)=(-1)^n\cdot C^{(3)}_n \cdot \frac{\binom{2n+j}{j}\cdot\binom{2n+1}{n+j+1}}{2n+1}\cdot (2n^2+n+1-j(n-1))\text{.}
\end{equation}
\end{lemma}

Let $Q(2n,m,a)$ denote the sum $\sum_{k=0}^{2n}(-1)^k\binom{2n}{k}^m\binom{a+k}{a}\binom{a+2n-k}{a}$.
It is known that $Q(2n,m,a)$ is divisible by $lcm(\binom{a+n}{a},\binom{2n}{n})$ for any non-negative integer $n$ and for any natural number $m$. See \cite [Thm.~12, p.\ 55]{JM3} and \cite[Remark 4, p.\ 860]{JM1}.
Furthermore, by using Eq.~(\ref{eq:9}), the following equations can be seen to hold:
\begin{align}
    M_{Q}(2n,j,0;a)&=(-1)^n\binom{a+n}{a}\binom{a+j}{j}\binom{a}{n-j}\label{eq:11}\text{,}\\
     M_{Q}(2n,j,1;a)&=(-1)^n\binom{2n}{n}\sum_{u=0}^{n-j}\binom{n}{j+u}\binom{j+u}{u}\binom{a+j+u}{j+u}\binom{a+n}{2n-j-u}\label{eq:12}\text{.}
\end{align}

For various applications of $M$ sums, see for example \cite[Background]{JM1}, \cite{JM3}, and \cite{JM4}. Recently, $M$ sums have been used for proving a new theorem  \cite[Thm.~\ 1, p.\ 5]{JM2} in number theory.

\section{Auxiliary Results}\label{s:3}

We shall also prove the following auxiliary propositions.

\begin{proposition}\label{p:1}
$T(n,j)$ is an integer for any non-negative integers $n$ and $j$ such that $j\leq n$. Furthermore, if $n$ is a natural number, then 
\begin{equation}\notag
T(n,j)=\frac{1}{n}\binom{n}{j}\binom{2n+j}{n-1}\text{.}
\end{equation}
\end{proposition}

\begin{proposition}\label{p:2}
Let $n$ be a non-negative integer. Then $\frac{2\binom{3n}{n}}{(n+1)(2n+1)}$ is an integer.
\end{proposition}

\begin{proposition}\label{p:3}
Let $n$ and $t$ be non-negative integers. Then $t\cdot \binom{2n+t}{t}$ is divisible by $2n+1$. 
\end{proposition}

\begin{proposition}\label{p:4}
Let $n$ and $t$ be non-negative integers. Then $\binom{3n}{n+t}\cdot \binom{2n+t}{2n}$ is divisible by $2n+1$.    
\end{proposition}

\begin{proposition}\label{p:5}
Let $n$ and $t$ be non-negative integers. Then $\binom{3n+1}{n+t+1}\cdot \binom{2n+t}{2n}$ is divisible by $2n+1$.    
\end{proposition}

The rest of the paper is structured as follows.
In Section $4$, we give a proof of Lemma \ref{l:1}. In Section $5$, we give proofs of Thms.~\ref{t:1} and \ref{t:2} and the proof of Proposition \ref{p:1}. In Section $6$, we prove Theorem \ref{t:3}, as well as Propositions \ref{p:2}, \ref{p:3}, \ref{p:4}, and \ref{p:5}. Finally, in Section $7$, we give  proofs of Thms.~\ref{t:5} and \ref{t:4}.

\section{Proof of Lemma \ref{l:1}}\label{s:4}
Let $S(2n,m)$ denote the sum from Eq.~(\ref{eq:3}). Obviously, $S(2n,m)$ corresponds to Definition \ref{def:1}, where
$F(2n,k)=(-1)^k\cdot C(2n,k)\cdot C(2n,2n-k)$.

By setting $t=0$ in Eq.~(\ref{eq:6}), we obtain that:
\begin{equation}\label{eq:13}
 M_S(2n,j,0)=\binom{2n-j}{j}\sum_{v=0}^{2n-2j}\binom{2n-2j}{v} F(2n,j+v)\text{,}\
 \end{equation}
 where $j$ and $t$ are non-negative integers such that $j\leq n$.

 By Eq.~(\ref{eq:1}), we can rewrite $F(2n,k)$ as
 \begin{equation}\label{eq:14}
   F(2n,k)= \frac{1}{(2n+1)^2}\cdot F_1(2n,k,1)\text{,}
 \end{equation}
 where $F_1(2n,k,1)$ is equal to 
 \begin{equation}\label{eq:15}
 (-1)^k \cdot (2n-k+1)(k+1) \cdot \binom{2n+k}{2n}\cdot\binom{4n-k}{2n}\text{.}
 \end{equation}
 
 Now, let $S_1(2n,m,1)$ denote the following sum $\sum_{k=0}^{2n}\binom{2n}{k}^m F_1(2n,k,1)$.

 By Eqs.~(\ref{eq:13}) and (\ref{eq:14}), we obtain that
 \begin{equation}\label{eq:16}
 M_S(2n,j,0)=\frac{1}{(2n+1)^2}\cdot M_{S_1}(2n,j,0;1)\text{.}
 \end{equation}

 More precisely, by Eq.~(\ref{eq:15}), the sum $S_1(2n,m,1)$ can be written as:
 \begin{equation}\label{eq:17}
     S_1(2n,m,1)=\sum_{k=0}^{2n}\binom{2n}{k}^m \binom{1+k}{1}\cdot \binom{1+2n-k}{1} \cdot (-1)^k \cdot \binom{2n+k}{2n}\cdot\binom{4n-k}{2n}\text{.}
 \end{equation}

 By setting $P:=S_1$, $a=1$, and $F_2(2n,k,1):=(-1)^k\cdot \binom{2n+k}{2n}\cdot\binom{4n-k}{2n}$ in Eq.~(\ref{eq:9}), we obtain that:
 \begin{equation}\label{eq:18}
M_{S_1}(2n,j,0;1)=\binom{1+j}{1}\sum_{l=0}^{1}\binom{2n-j+l}{l}\binom{2n-j}{1-l}M_{S_2}(2n,j+1-l,0;1)\text{,}
\end{equation}
where 
\begin{equation}\label{eq:19}
S_2(2n,m,1)=\sum_{k=0}^{2n}\binom{2n}{k}^m (-1)^k\cdot \binom{2n+k}{2n}\cdot\binom{4n-k}{2n}\text{.}
\end{equation}

Note that: 
\begin{equation}\label{eq:20}
 S_2(2n,m,1)=Q(2n,m,2n)\text{.}   
\end{equation}

Hence, by Eq.~(\ref{eq:5}), it follows
\begin{equation}\label{eq:21}
   M_{S_2}(2n,j,0;1)=M_Q(2n,j,0;2n)\text{.}
\end{equation}

Therefore, by Eq.~(\ref{eq:18}), it follows that $M_{S_1}(2n,j,0;1)$ is equal to

\begin{equation}\label{eq:22}
    (1+j)\sum_{l=0}^{1}\binom{2n-j+l}{l}\binom{2n-j}{1-l}M_{Q}(2n,j+1-l,0;2n)\text{.}
\end{equation}

The first summand in Eq.~(\ref{eq:22}), for $l=0$, is equal to
\begin{equation}\label{eq:23}
  (1+j)\cdot (2n-j)\cdot M_{Q}(2n,j+1,0;2n)\text{.}
\end{equation}

By setting $a:=2n$ and $j:=j+1$ in Eq.~(\ref{eq:11}), it follows that
\begin{equation}
M_{Q}(2n,j+1,0;2n)=(-1)^n\binom{3n}{2n}\binom{2n+j+1}{j+1}\binom{2n}{n-j-1}\text{.}\label{eq:24}
\end{equation}

Similarly, the second summand in Eq.~(\ref{eq:22}), for $l=1$, is equal to
\begin{equation}\label{eq:25}
  (1+j)\cdot (2n-j+1)\cdot M_{Q}(2n,j,0;2n)\text{.}
\end{equation}

By setting $a:=2n$ in Eq.~(\ref{eq:11}), we obtain that
\begin{equation}
M_{Q}(2n,j,0;2n)=(-1)^n\binom{3n}{2n}\binom{2n+j}{j}\binom{2n}{n-j}\text{.}\label{eq:26}
\end{equation}

Finally, by Eqs.~(22),(23), (24), and (25), we obtain that $M_{S_1}(2n,j,0;1)$ is equal to
\begin{equation}\label{eq:27}
    (1+j)(-1)^n\binom{3n}{2n}\left ( (2n-j)\binom{2n+j+1}{j+1}\binom{2n}{n-j-1}+(2n-j+1)\binom{2n+j}{j}\binom{2n}{n-j}\right )\text{.}
\end{equation}

We shall use two well-known binomial formulae:
\begin{align}
    \binom{n+j+1}{j+1}&=\frac{n+j+1}{j+1}\binom{n+j}{j}\text{,}\label{eq:28}\\
    \binom{2n}{n-j-1}&=\frac{n-j}{n+1+j}\cdot \binom{2n}{n-j}\text{.}\label{eq:29}
\end{align}

We have that $(2n-j)\binom{2n+j+1}{j+1}\binom{2n}{n-j-1}$ is equal to:

\begin{equation}\label{eq:30}
    (2n-j)\cdot \frac{2n+j+1}{j+1}\cdot \frac{n-j}{n+1+j}\cdot \binom{2n+j}{j}\cdot \binom{2n}{n-j}\text{.}
\end{equation}

After cancellation of $1+j$, it follows that $M_{S_1}(2n,j,0;1)$ is equal to 
\begin{equation}\label{eq:31}
    (-1)^n\binom{3n}{2n}\binom{2n+j}{j}\cdot \binom{2n}{n-j}\bigl{(}(2n-j)\cdot (2n+j+1)\cdot \frac{n-j}{n+1+j}+(2n-j+1)(j+1)\bigr{)}\text{.}
\end{equation}

It can be shown, after some computations, that 
$(2n-j)\cdot (2n+j+1)\cdot \frac{n-j}{n+1+j}+(2n-j+1)(j+1)$ is equal to
\begin{equation}\label{eq:32}
\frac{4n^3+2n^2(2-j)+n(j+3)+(j+1)}{n+j+1}\text{.}
\end{equation}

Therefore, $M_{S_1}(2n,j,0;1)$ is equal to

\begin{equation}\label{eq:33}
    (-1)^n\binom{3n}{2n}\binom{2n+j}{j}\cdot \binom{2n}{n-j}\cdot \frac{4n^3+2n^2(2-j)+n(j+3)+(j+1)}{n+j+1}\text{.}
\end{equation}

By Eq.~(\ref{eq:16}), it follows that $M_S(2n,j,0)$ is equal to
\begin{equation}\label{eq:34}
\frac{1}{(2n+1)^2}\cdot  (-1)^n\binom{3n}{2n}\binom{2n+j}{j}\cdot \binom{2n}{n-j}\cdot \frac{4n^3+2n^2(2-j)+n(j+3)+(j+1)}{n+j+1}\text{.}
\end{equation}

Note that polynomial $4n^3+2n^2(2-j)+n(j+3)+(j+1)$ is divisible by $2n+1$, since
\begin{equation}\label{eq:35}
4n^3+2n^2(2-j)+n(j+3)+(j+1)=(2n^2+n(1-j)+(j+1))\cdot (2n+1)\text{.}
\end{equation}

By Eq.~(\ref{eq:35}), the sum $M_S(2n,j,0)$ from Eq.~(\ref{eq:34}) becomes
\begin{equation}\label{eq:36}
(-1)^n\frac{1}{2n+1}\cdot \binom{3n}{2n}\binom{2n+j}{j}\cdot \binom{2n}{n+j}\cdot \frac{1}{n+j+1}\cdot (2n^2+n(1-j)+(j+1) \text{.}
\end{equation}
It is readily verified that \cite[Eq.~(1.1), p.\ 5]{TC}:
\begin{equation}\label{eq:37}
\binom{2n}{n+j}\cdot \frac{1}{n+j+1}=\frac{1}{2n+1} \cdot \binom{2n+1}{n+j+1}\text{.}
\end{equation}

By using Eq.~(\ref{eq:37}), the sum $M_S(2n,j,0)$ from Eq.~(\ref{eq:36}) becomes
\begin{equation}\label{eq:38}
(-1)^n\cdot\frac{1}{2n+1}\cdot \binom{3n}{2n} \cdot \frac{\binom{2n+j}{j}\cdot\binom{2n+1}{n+j+1}}{2n+1}\cdot (2n^2+n+1-j(n-1) \text{.}
\end{equation}

By setting $C^{(3)}_n=\frac{1}{2n+1}\cdot \binom{3n}{2n}$, we obtain that
\begin{equation}\label{eq:39}
M_S(2n,j,0)=(-1)^n\cdot C^{(3)}_n \cdot \frac{\binom{2n+j}{j}\cdot\binom{2n+1}{n+j+1}}{2n+1}\cdot (2n^2+n+1-j(n-1))\text{.}
\end{equation}

The last equation above proves Eq.~(\ref{eq:10}). This completes the proof of Lemma \ref{l:1}.

\section{Proofs of Theorem \ref{t:1}, Proposition \ref{p:1}, and Theorem \ref{t:2}}\label{s:5}

We give a proof of Theorem \ref{t:1}.

Let us recall that $T(n,j)$ is equal to $\frac{\binom{2n+j}{j}\cdot\binom{2n+1}{n+j+1}}{2n+1}$.

\subsection{Proof of Theorem \ref{t:1}}\label{s:5.1}
Obviously, $T(n,0)=C_n$. By setting $j:=0$ in Eq.~(\ref{eq:10}), we obtain that
\begin{equation}\label{eq:40}
 M_S(2n,0,0)=(-1)^n\cdot C^{(3)}_n \cdot C_n \cdot (2n^2+n+1)\text{.}   
\end{equation}

By setting $m:=1$ in Eq.~(\ref{eq:7}), it follows that:
\begin{equation}\label{eq:41}
    S(2n,1)=M_S(2n,0,0)\text{.}
\end{equation}

Finally, by Eqs.~(\ref{eq:40}) and (\ref{eq:41}), it follows that
\begin{equation}\notag
S(2n,1)=(-1)^n\cdot C^{(3)}_n \cdot C_n \cdot (2n^2+n+1)\text{.}
\end{equation}

The last equation above is exactly Eq.~(\ref{eq:4}). This completes the proof of Theorem \ref{t:1}.

Next, we give a proof of Proposition \ref{p:1}.

\subsection{Proof of Proposition \ref{p:1}}\label{s:5.2}

Clearly, $T(0,0)=1$.
Now, let us assume that $n$ is a natural number.

It is readily verified \cite[Eq.~(1.1),p.\ 5]{TC} that
\begin{equation}\label{eq:42}
    \binom{2n+1}{n+j+1}=\frac{2n+1}{n+j+1}\binom{2n}{n+j}\text{.}
\end{equation}
We have
\begin{align}
    T(n,j)&=\binom{2n+j}{2n}\cdot \binom{2n}{n+j}\cdot \frac{1}{n+j+1}&&(\text{by using Eq.~(\ref{eq:42})})\notag\\
    &=\binom{2n+j}{n+j}\cdot \binom{n}{j} \cdot \frac{1}{n+j+1}&& (\text{by using \cite[Eq.~(1.4), p.\ 5]{TC}})\notag\\
    &=\binom{2n+j}{n+j}\cdot \frac{1}{n+j+1}\cdot \binom{n}{j}\notag \\
    &=\frac{1}{n}\cdot \binom{2n+j}{n+j+1}\cdot \binom{n}{j}&& (\text{since } n \text{ is a natural number})\notag\\
    &=\frac{1}{n} \cdot \binom{2n+j}{n-1}\cdot \binom{n}{j}\text{.}\label{eq:43}
\end{align}

By Eq.~(\ref{eq:43}), it follows that

\begin{equation}\label{eq:44}
    n\cdot\binom{2n+j}{j}\cdot\binom{2n+1}{n+j+1}=(2n+1) \cdot \binom{2n+j}{n-1}\cdot \binom{n}{j}\text{.}
\end{equation}

Note that the integers $n$ and $2n+1$ are relatively prime. Therefore, it follows from Eq.~(\ref{eq:44}) that 
$n$ must divide $\binom{2n+j}{n-1}\cdot \binom{n}{j}$, and $2n+1$ must divide $\binom{2n+j}{j}\cdot\binom{2n+1}{n+j+1}$.
This completes the proof of Proposition \ref{p:1}.

Now, we are ready to prove Theorem \ref{t:2}.

\subsection{Proof of Theorem \ref{t:2}}\label{s:5.3}

By using Proposition \ref{p:1}, we know that $T(n,j)$ is an integer
for any non-negative integer $n$ and any $j$ such that $j\leq n$.
By Eq.~(\ref{eq:10}), $C^{(3)}_n=\frac{1}{2n+1}\cdot \binom{3n}{2n}$ divides $M_S(2n,j,0)$
for any non-negative $j$ such that $j\leq n$. By Eq.~(\ref{eq:7}), $C^{(3)}_n$ divides the sum $M_S(2n,j,t)$
for any non-negative integer $t$. By Eq.~(\ref{eq:6}), the sum $S(2n,m)$ is divisible by $C^{(3)}_n$ for any non-negative integer $n$, and for any natural number $m$. This completes the proof of Theorem \ref{t:2}.

\section{Proof of Theorem \ref{t:3}}\label{s:6}
Our proof of Theorem \ref{t:3} consists of two parts. In the first part, we prove that the sum $S(2n,1)$ is divisible by $\binom{2n}{n}$ for any non-negative integer $n$. In the second part, we prove that $S(2n,m)$ is divisible by $\binom{2n}{n}$ for any non-negative integer $n$ and for any natural number $m$ such that $m\geq 2$.

We begin with the proof of Proposition \ref{p:2}.
\subsection{Proof of Proposition \ref{p:2}}\label{s:6.1}
Let us prove that $\frac{2\binom{3n}{n}}{n+1}$ is an integer.

  We have
  \begin{align}
\binom{3n}{n+1}\cdot (n+1)&=\binom{3n}{n+1}\cdot \binom{n+1}{n}\text{,}\notag\\
&=\binom{3n}{n}\cdot \binom{2n}{1}\text{,}&&(\text{By Newton's identity } \cite[Eq.~(1.4),p.\ 5]{TC})\notag\\
&=\binom{3n}{n}\cdot (2n)\text{.}\label{eq:45}
  \end{align}
Since $n+1$ and $n$ are relatively prime integers, it follows from Eq.~(\ref{eq:45}),
that $n+1$ must divide $2\binom{3n}{n}$.

Next, from the definition of Fuss-Catalan number $C^{(3)}_n$ of order three, we know that $2n+1$ must divide the binomial coefficient $\binom{3n}{n}$.
Obviously, then $2n+1$ must also divide $2\binom{3n}{n}$.

Finally, due to the fact that $n+1$ and $2n+1$ are relatively prime, it follows that $(n+1)(2n+1)$ must divide $2\binom{3n}{n}$.
Therefore, $\frac{2\binom{3n}{n}}{(n+1)(2n+1)}$ must be an integer. This completes the proof of Proposition \ref{p:2}.

\subsection{Proof of the First Part of Theorem \ref{t:3}}\label{s:6.2}

By Eq.~(\ref{eq:4}), it follows that
\begin{equation}\label{eq:46}
S(2n,1)=(-1)^n\cdot \frac{2n^2+n+1}{(n+1)(2n+1)}\cdot \binom{3n}{n}\cdot \binom{2n}{n}\text{.}
\end{equation}

Let us prove that $\frac{2n^2+n+1}{(n+1)(2n+1)}\cdot \binom{3n}{n}$ is an integer for any non-negative integer $n$.

It is readily verified that

\begin{align}
   \frac{2n^2+n+1}{(n+1)(2n+1)}\cdot \binom{3n}{n}&=\frac{2n^2}{(n+1)(2n+1)}\cdot \binom{3n}{n}+C^{(3)}_n\text{,}\notag\\
   &=n^2\cdot \frac{2}{(n+1)(2n+1)}\cdot \binom{3n}{n}+C^{(3)}_n\text{.}\label{eq:47}
\end{align}

By Proposition \ref{p:2},  $\frac{2n^2}{(n+1)(2n+1)}\cdot \binom{3n}{n}$ must be an integer. By using Eq.~(\ref{eq:47}) and the fact that the Fuss-Catalan number of order three $C^{(3)}_n$ is an integer, it follows that $\frac{2n^2+n+1}{(n+1)(2n+1)}\cdot \binom{3n}{n}$ is an integer.

By using Eq.~(\ref{eq:47}) and the fact that $\frac{2n^2+n+1}{(n+1)(2n+1)}\cdot \binom{3n}{n}$ is an integer, it follows that $S(2n,1)$ is divisible by the central binomial coefficient $\binom{2n}{n}$. This completes the first part of the proof of Theorem \ref{t:3}.

Now, we give proofs of the remaining propositions.

\subsection{Proof of Proposition \ref{p:3}}\label{s:6.3}

For $t=0$, Proposition \ref{p:3} is true. Let us therefore assume that $t$ is a natural number. 

We have 
\begin{align}
    \binom{2n+t}{t}\cdot t&=\binom{2n+t}{t}\cdot \binom{t}{t-1}\text{,}\notag\\
    &=\binom{2n+t}{t-1}\cdot \binom{2n+1}{1}\text{.}&& (\text{By } \cite[Eq.~(1.4),p.\ 5]{TC})\label{eq:48}
\end{align}

By Eq.~(\ref{eq:48}), it follows that
\begin{equation}\label{eq:49}
 \binom{2n+t}{t}\cdot t=(2n+1)\cdot \binom{2n+t}{t-1}\text{.}    
\end{equation}

Finally, Eq.~(\ref{eq:49}) completes the proof of Proposition \ref{p:3}.

\subsection{Proof of Proposition \ref{p:4}}\label{s:6.4}

Let us consider the integer $(3n+1)\binom{3n}{n+t}\cdot \binom{2n+t}{2n}$.

We have

\begin{align}
 (3n+1)\binom{3n}{n+t}\cdot \binom{2n+t}{2n}&=\binom{3n+1}{3n}\cdot \binom{3n}{n+t}\cdot \binom{2n+t}{2n}\text{,}\notag\\
 &=\binom{3n+1}{n+t}\cdot \binom{2n+1-t}{2n-t}\cdot \binom{2n+t}{2n}\text{,}\notag\\
 &=\binom{3n+1}{n+t}\cdot (2n+1-t)\cdot \binom{2n+t}{2n}\text{.}\label{eq:50}
\end{align}

By Eq.~(\ref{eq:50}) and the symmetry of binomial coefficients, it follows that
\begin{equation}\label{eq:51}
(3n+1)\binom{3n}{n+t}\cdot \binom{2n+t}{2n}=(2n+1)\cdot\binom{3n+1}{n+t}\cdot \binom{2n+t}{2n}-\binom{3n+1}{n+t}\cdot t \cdot \binom{2n+t}{t}\text{.}
\end{equation}

By using Eq.~(49), Eq.~(\ref{eq:51}) becomes 
\begin{equation}
(3n+1)\binom{3n}{n+t}\cdot \binom{2n+t}{2n}=(2n+1)\cdot \binom{3n+1}{n+t} \left(\binom{2n+t}{2n}-\binom{2n+t}{t-1}\right)\text{.}\label{eq:52}
\end{equation}

Due to the fact that $2n+1$ and $3n+1$ are relatively prime, it follows from Eq.~(\ref{eq:52}) that $2n+1$ must divide $\binom{3n}{n+t}\cdot \binom{2n+t}{2n}$. This completes the proof of Proposition \ref{p:4}.

\subsection{Proof of Proposition \ref{p:5}}\label{s:6.5}

Let us consider the integer $(3n+2)\binom{3n+1}{n+t+1}\cdot \binom{2n+t}{2n}$.

We shall use the well-known formula \cite[Eq.~(1.4),p.\ 5]{TC}. We have

\begin{align}
 (3n+2)\binom{3n+1}{n+t+1}\cdot \binom{2n+t}{2n}&=\binom{3n+2}{3n+1}\cdot \binom{3n+1}{n+t+1}\cdot \binom{2n+t}{2n}\text{,}\notag\\
 &=\binom{3n+2}{n+t+1}\cdot \binom{2n+1-t}{2n-t}\cdot \binom{2n+t}{2n}\text{,}\notag\\
 &=\binom{3n+2}{n+t+1}\cdot (2n+1-t)\cdot \binom{2n+t}{2n}\text{.}\label{eq:53}
\end{align}

By using Eq.~(\ref{eq:49}) and the symmetry of binomial coefficients, Eq.~(\ref{eq:53}) becomes
\begin{equation}
(3n+2)\binom{3n+1}{n+t+1}\cdot \binom{2n+t}{2n}=(2n+1)\cdot \binom{3n+2}{n+t+1} \left(\binom{2n+t}{2n}-\binom{2n+t}{t-1}\right) \text{.}\label{eq:54}
\end{equation}

Due to the fact that integers $2n+1$ and $3n+2$ are relatively prime, it follows from Eq.~(\ref{eq:54}) that $2n+1$ must divide
$\binom{3n+1}{n+t+1}\cdot \binom{2n+t}{2n}$. This completes the proof of Proposition \ref{p:5}.

\begin{remark}\label{r:1}
By using Propositions \ref{p:4} and \ref{p:5}, it also follows that
$\binom{3n}{n+t+1}\cdot \binom{2n+t}{2n}$ is divisible by $2n+1$. One can use Pascal's formula \cite[Theorem 1.1, p.\ 5]{TC}:
\begin{equation}\label{eq:55}
\binom{3n+1}{n+t+1}=\binom{3n}{n+t+1}+\binom{3n}{n+t}\text{.}
\end{equation}
\end{remark}

\subsection{Proof of the Second Part of Theorem \ref{t:3}}\label{s:6.6}

We shall prove that the sum $M_S(2n,j,1)$ is divisible by $\binom{2n}{n}$ for any non-negative integers $n$ and any $j$ such that $j\leq n$,
where $S(2n,m)$ is our main sum from Eq.~(\ref{eq:3}).

By setting $t:=0$ in Eq.~(\ref{eq:8}), it follows that
\begin{equation}
M_S(2n,j,1)=\binom{2n}{j}\sum_{u=0}^{n-j}\binom{2n-j}{u}\cdot M_S(2n,j+u,0)\text{.}\label{eq:56}   
\end{equation}

By setting $j:=j+u$ in Eq.~(\ref{eq:10}), it follows that $M_S(2n,j+u,0)$ is equal to
\begin{equation}\label{eq:57}
(-1)^n\binom{3n}{n} \cdot \frac{1}{2n+1}\cdot \frac{\binom{2n+j+u}{j}\cdot\binom{2n+1}{n+j+u+1}}{2n+1}\cdot (2n^2+n+1-(j+u)(n-1))\text{.}    
\end{equation}

By using Eq.~(\ref{eq:42}), we obtain that
\begin{equation}\label{eq:58}
\frac{\binom{2n+1}{n+j+u+1}}{2n+1}=\frac{\binom{2n}{n+j+u}}{n+j+u+1}\text{.}
\end{equation}

By using Eqs.~(\ref{eq:57}) and (\ref{eq:58}), it follows that the summand $\binom{2n}{j}\binom{2n-j}{u}M_S(2n,j+u,0)$ from Eq.~(\ref{eq:56}) is equal to
\begin{equation}\label{eq:59}
\frac{(-1)^n}{2n+1}\binom{3n}{2n}\binom{2n}{j}\binom{2n-j}{u}\frac{1}{n+j+u+1}\binom{2n+j+u}{j+u}\binom{2n}{n+j+u}\cdot (2n^2+n+1-(j+u)(n-1))\text{.}
\end{equation}

By using \cite[Eq.~(1.4),p.\ 5]{TC}, it is readily verified that
\begin{equation}\label{eq:60}
\binom{3n}{2n}\binom{2n}{j}\binom{2n-j}{u}=\binom{3n}{n+j+u}\binom{n+j+u}{n}\binom{j+u}{u}\text{.}  
\end{equation}

Similarly, by using \cite[Eq.~(1.4),p.\ 5]{TC}, it follows that
\begin{equation}\label{eq:61}
    \binom{2n}{n+j+u}\binom{n+j+u}{n}=\binom{2n}{n}\cdot \binom{n}{j+u}\text{.}
\end{equation}

By using Eqs.~(\ref{eq:60}), (\ref{eq:61}), and the substitution $t=j+u$, the inner summand of Eq.~(\ref{eq:59}) becomes:
\begin{equation}\label{eq:62}
 \frac{(-1)^n}{2n+1}\binom{3n}{n+t}\binom{2n+t}{t}\cdot\binom{2n}{n}\binom{n}{t}\binom{t}{j} \frac{1}{n+t+1}(2n^2+n+1-t(n-1)\text{.}   
\end{equation}

Let us prove that the number 
\begin{equation}\label{eq:63}
N(n,j,t)=\frac{(-1)^n}{2n+1}\binom{3n}{n+t}\binom{2n+t}{t}\binom{n}{t}\binom{t}{j} \frac{1}{n+t+1}(2n^2+n+1-t(n-1)
\end{equation}
is an integer for any non-negative integers $n$, $j$, and $t$ such that $j\leq t \leq n$.

Note that 
\begin{equation}\label{eq:64}
\frac{2n^2+n+1-t(n-1)}{n+t+1}=\frac{n(2n-t)}{n+t+1}+1\text{.}
\end{equation}

It is readily verified that
\begin{equation}\label{eq:65}
\binom{3n}{n+t}\cdot \frac{2n-t}{n+t+1}=\binom{3n}{n+t+1}\text{.}
\end{equation}

By using Eqs.~(\ref{eq:64}) and (\ref{eq:65}), it follows that
\begin{equation}
N(n,j,t)=N_1(n,j,t)+N_2(n,j,t)\text{,}\label{eq:66}    
\end{equation}
where 
\begin{align}
 N_1(n,j,t)&=(-1)^n\frac{\binom{3n}{n+t}\binom{2n+t}{2n}}{2n+1}\binom{n}{t}\binom{t}{j} \text{,}\label{eq:67}\\
 N_2(n,j,t)&=n(-1)^n\frac{\binom{3n}{n+t+1}\binom{2n+t}{2n}}{2n+1}\binom{n}{t}\binom{t}{j} \text{.}\label{eq:68}
\end{align}

By Proposition \ref{p:4}, $N_1(n,j,t)$ is an integer.
By Remark \ref{r:1}, $N_2(n,j,t)$ is also an integer.

Finally, by Eq.~(\ref{eq:66}), it follows that $N(n,j,t)$ is also an integer.

Furthermore, by using Eqs.~(\ref{eq:56}) and (\ref{eq:63}), it follows that
\begin{equation}\label{eq:69}
M_S(2n,j,1)=\binom{2n}{n}\sum_{u=0}^{n-j}N(n,j,t)\text{,}   
\end{equation}
where $t=j+u$.

Since $N(n,j,t)$ is an integer, it follows from Eq.~(\ref{eq:69}) that $M_S(2n,j,1)$ is divisible by the central binomial coefficient $\binom{2n}{n}$ for any non-negative integers $n$ and $j$ such that $j\leq n$.

By using Eq.~(\ref{eq:8}) and the induction principle, it can be shown that the sum $M_S(2n,j,t)$ is divisible by $\binom{2n}{n}$
for any non-negative integers $n$ and $j$  such that $j\leq n$, and for any natural number $t$.

By using Eq.~(\ref{eq:7}), it follows that the sum $S(2n,m)$ is divisible by $\binom{2n}{n}$ for any non-negative integer $n$ and for any natural number $m$ such that $m\geq 2$.
This proves the second part of the proof of Theorem \ref{t:3}.

\section{Proof of Theorem \ref{t:5}}\label{s:7}

We shall use the induction principle on the pairs $(n,m)$, where $m \geq 2n$. We say that $(n_1,m_1)<(n_2,m_2)$ if $n_1<n_2$ or $ n_1=n_2$ and $m_1<m_2$.

For $n=1$ and for any natural number $m\geq 2$, it is readily verified that $Schr(1,m,j,2)=m-1$, where $j=0$ or $j=1$. Therefore, Eq.~(\ref{eq:5}) is true for all pairs $(1,m)$ where $m\geq 2$. Note that $Schr(n,m,j,2)=0$ if $m <2n$ and $Schr(n,m,j,2)=0$ if $j>n$ or $j<0$. Thus Eq.~(\ref{eq:5}) also holds whenever $m=2n-1$.

Let us suppose that $n$ is a fixed natural number greater than $1$, and let $m$ be an arbitrary natural number such that $m\geq 2n$.

Now, we use a substitution $m=2n+p$, where $p$ is a non-negative integer. We have that $Schr(n,2n+p,j,2)$ is equal to
\begin{equation}\label{eq:70}
Schr(n-1,2n+p,j-1,2)+Schr(n-1,2n+p-1,j,2)+Schr(n,2n+p-1,j,2)\text{,}
\end{equation}
where $0\leq j\leq n$.

Let us suppose that Eq.~(\ref{eq:5}) is true for all pairs $(n_1, m_1)$ such that $(n_1, m_1)<(n,2n+p)$.

Therefore, by using the induction hypothesis, we know that:
\begin{align}
Schr(n-1,2n+p,j-1,2)&=\frac{p+3}{n-1}\binom{n-1}{j-1}\binom{2n+p+j-1}{n-2}\text{,}\label{eq:71}\\
Schr(n-1,2n+p-1,j,2)&=\frac{p+2}{n-1}\binom{n-1}{j}\binom{2n+p+j-1}{n-2}\text{,}\label{eq:72}\\
Schr(n,2n+p-1,j,2)&=\frac{p}{n}\binom{n}{j}\binom{2n+p+j-1}{n-1}\text{.}\label{eq:73}
\end{align}

It can be shown that by adding Eqs.~(\ref{eq:71}) and (\ref{eq:72}), we obtain that 
\begin{equation}\label{eq:74}
Schr(n-1,2n+p,j-1,2)+Schr(n-1,2n+p-1,j,2)=\frac{1}{n}\binom{n}{j}\binom{2n+p+j-1}{n-2}\frac{2n+pn+j}{n-1}\text{.}
\end{equation}

Furthermore, by adding Eqs.~(\ref{eq:74}) and (\ref{eq:73}), 
we obtain that Eq.~(\ref{eq:70}) becomes
\begin{equation}\label{eq:75}
 \frac{1}{n}\binom{n}{j}\binom{2n+p+j-1}{n-2}\frac{2n+pn+j}{n-1}+ \frac{p}{n}\binom{n}{j}\binom{2n+p+j-1}{n-1}\text{.}  
\end{equation}

After some calculations, it can be shown that Eq.~(\ref{eq:75}) becomes
\begin{equation}\label{eq:76}
\frac{p+1}{n}\binom{n}{j}\binom{2n+p+j}{n-1}\text{.}
\end{equation}

By using the induction principle, it follows that Eq.~(\ref{eq:5}) is true for all natural numbers $n$ and $m$ such that $m\geq 2n$.
This completes the proof of Theorem \ref{t:5}.

\begin{remark}\label{r:2}
Theorem \ref{t:4} follows from Theorem \ref{t:5} by setting $m=2n$ in Eq.~(\ref{eq:5}) and by using Proposition \ref{p:1}.
\end{remark}

\begin{remark}\label{r:3}
Let $n$, $m$, and $l$ be natural numbers such that $m\geq ln$.
Let the number $Schr(n,m,j,l)$ denote the number of all lattice paths
in the plane from the point $(0,0)$ to the point $(n,m)$, using steps $(1,0)$, $(0,1)$, and $(1,1)$,
that never go below the line $y=lx$ with exactly $j$ steps of the form $(1,0)$. By a using similar idea to the one in the proof of Theorem \ref{t:5},
it can be shown that 
\begin{equation}\label{eq:79}
Schr(n,m,j,l)=\frac{m-ln+1}{n}\binom{n}{j}\binom{m+j}{n-1}\text{.}
\end{equation}
Furthermore, let the number $Schr(n,m,l)$ denote the number of all lattice paths
in the plane from the point $(0,0)$ to the point $(n,m)$, using steps $(1,0)$, $(0,1)$, and $(1,1)$,
that never go below the line $y=lx$.
By using a combinatorial interpretation for the numbers $Schr(n,m,l,j)$ in the Eq.~(\ref{eq:79}), we obtain a new combinatorial proof of the following known \cite[Theorem 2.9, p.\ 6]{JS} formula
\begin{equation}\label{eq:80}
Schr(n,m,l)=\frac{m-ln+1}{n}\sum_{j=0}^{n}\binom{n}{j}\binom{m+j}{n-1}\text{.}
\end{equation}    
\end{remark}

\section*{Acknowledgments}
I want to thank Professor Tomislav Do\v{s}li\'{c} for inviting me to the second and third Croatian Combinatorial Days in Zagreb. Also, I am very grateful to Professor David Dol\v{z}an for his valuable comments which improved the quality of this manuscript.

\bigskip
\hrule
\bigskip

\noindent 2010 {\it Mathematics Subject Classification}:
Primary  05A10 ; Secondary 05A19.

\noindent\emph{Keywords:} Catalan Triangle, Fuss-Catalan Number of Order Three, Catalan Number, Central Binomial Coefficient, $M$ sum, Schr\"{o}der Number of Order Two, Lattice Path, Induction Principle.

\bigskip
\hrule 
\bigskip

\noindent

\end{document}